**Short Paper**

# Evaluating the Effect of Activity-Based Method of Teaching Mathematics on Nigerian Secondary School Students' Achievement in Mathematics

Chinelo Blessing Oribhabor
Department of Guidance and Counseling, Faculty of Arts and Education, University of Africa, Toru-Orua, Bayelsa State, Nigeria
chiblessing42004@yahoo.co.uk; chinelo.oribhabor@uat.edu.ng

**Abstract**

   Mathematics is a compulsory subject in Nigerian secondary schools, and the subject plays an important role in the scientific and technological growth and development of the nation. A shortfall in the knowledge of the students in Mathematics means that the goal may not be realized, hence the need to improve teaching methods for solving the problem of poor performance in the subject. This study evaluated the effect of the activity-based teaching method on the students' achievement in secondary school Mathematics. The design of the study was a quasi-experimental pretest-posttest research design using intact classes. A sample of 96 students from two senior secondary schools in Calabar metropolis, Cross Rivers State, Nigeria was used in the study. The instrument for data collection was a 20-item Essay Mathematics Test (MT) developed by the researcher to measure students' achievement in Mathematics. The reliability coefficient for the instrument using the Cronbach coefficient alpha was 0.86. Mean and the standard deviation was used to answer the research question while Analysis of Covariance (ANCOVA) was used to explore the effect of the teaching method at a 0.05 level of significance. Finding revealed that there was a significant difference in the Mathematics performance between the posttest mean scores of the students who were exposed to activity-based teaching methods (experimental) and those that were taught with lecture method (control) groups after controlling for the effect of the pre-test on Mathematics scores. The result indicated that students taught using the activity method performed better than those taught using the lecture method. This study implies that students will benefit greatly from teaching outcomes in schools if the activity method is used as a pedagogical approach in Mathematics instructional delivery. The paper recommends among others that secondary school Mathematics teachers should be trained and re-trained to update their knowledge in the use of activity-based teaching for making the teaching and learning of Mathematics more interesting and rewarding.

## INTRODUCTION

Mathematics is used either consciously or unconsciously in various aspects of life and is the basic tool for industrialization and national development. It is also recognized to play a vital role in contemporary society, making it indispensable for the existence of any nation (Asante, 2010). The importance of Mathematics can be seen in its applications to science and technology, medicine, the economy, the environment, and in public decision making. The implication is that for effective functioning in society, there is a need for all the citizens to study and understand science and mathematics. The imperative role of mathematics notwithstanding, Mathematics is one of the most poorly taught, widely hated, and abysmally understood subjects in our schools (Ali, et al., 2010). This is evident in the persistent poor performance of Nigerian secondary school students in external examinations such as the Senior School Certificate Examination (SSCE).

Analysis of students' achievement in Mathematics shows that the performance of Nigerian students is not encouraging (Imoko & Jimin, 2016). According to the Naija News (2019), 43.30% of candidates that wrote the 2018 West African Examination Council (WAEC) in Cross River State, Nigeria obtained credits and above in Mathematics, this shows the abysmal performance in the state. WAEC Chief Examiner for the private candidates (2018) further observed that candidates were weak in the areas such as Mensuration, Translation of word problem to form equations, Sequence and Series, Inequality, Histogram, Circle Geometry, Matrices, and Construction. This has always led to poor academic achievement and has been a great thing of worry for all stakeholders, such as parents, teachers, educational psychologists, counselors, government, and society at large. This is because mathematics is today, part of the basic requirements for entrance into tertiary institutions. This consistent poor performance by students in mathematics calls for serious national action to remedy the situation.

Federal Ministry of Education, Nigeria, including the Cross River State Ministry of Education has made efforts to reduce the incessant poor performance of students in mathematics by embarking on several programs, such as the introduction of new mathematics textbooks, different mathematics competitions, workshops, and seminars, yet the problems seem to continue. Moreover, attempts to find a solution to the incessant failure rates of students have made researchers in mathematics education to consider several factors. Poor academic achievement in mathematics could be attributed to many factors among which are the teachers' teaching methods. This means that mathematics concepts cannot be well understood if students are not taught with an effective teaching method. Many researchers such as Kajuru and Popoola (2010), Adewumi (2010), Ifeanacho (2012), and Agommuoh and Ifeanacho (2013) have blamed this poor achievement in mathematics on the use of inappropriate teaching method which might lead to a lack of interest and retention of mathematical concepts. In the teaching of mathematics, strategies that involve critical thinking and the generation of innovative minds need to be employed. Unfortunately, the teaching of mathematics in Nigeria, and particularly in Cross River State is mainly theoretical and talk and chalk, with fewer hands-on activities and fewer methods that directly engage the student to capacity build the ability to think and solve problems independently.

It has also been reported by Mtsem (2011) that the teaching method affects the responses of students and determines whether they are interested, motivated, and involved in a lesson in such a way as to engage in good learning. What constitutes good teaching and learning of school subjects is the use of appropriate methods of teaching. Therefore, establishing a more effective teaching method in the teaching of Mathematics could offer a solution to poor performance in Mathematics. Saga (2015) asserts that if the existing strategy of teaching is not yielding results, then, other teaching strategies should be adopted. Given this, there is a need for a paradigm shift to using methods that are appropriate, relevant, and effective to the subject. Saribas (2009) cited in Mvula (2020) opined that when appropriate learning methods are used, the learners develop proper attitudes and skills during the learning process. Hence there is a need to change from the conventional way of teaching Mathematics to effective ones. The results of this study will provide useful insights to the Mathematics teachers, school management, and the ministries of secondary, Tertiary, and General Education respectively on the differential effectiveness of Mathematics teaching methods on students' academic performance in Mathematics. Also, academicians and scholars will use the study as a source of information regarding the effective learning of Mathematics. Finally, the study will add value to the literature on methods of enhancing the learning of Mathematics among secondary students.

This study is conducted to provide a framework analysis on the effect of activity-based teaching methods on the senior secondary school students' Mathematics achievement for the topics of 'Circle Geometry' and 'Construction' in comparison to the lecture method of teaching. The activity-based method was chosen, to help the secondary students to understand the two selected topics (Circle Geometry and Construction) very well.

Specifically, this study aims at finding answers to the following question:

1. Is there a difference between the Mathematics performance mean scores of students in the experimental (activity-based method) and control (lecture method) groups when the moderating effect of the pretest is controlled?

*Hypothesis*

The following hypothesis was tested at the 0.05 significance level:

1. There is no significant difference between the posttest Mathematics performance mean scores of students in the experimental (activity-based method) and control (lecture method) groups when the moderating effect of the pretest is controlled.

## LITERATURE REVIEW

There are several teaching methods but for this research, the researcher condensed them to the Lecture method and activity-based method. Activity-based learning (ABL) is an instructional approach that emphasizes students' active learning through various activities to develop the three domains of learning (cognitive, affective, and psychomotor) equally (Priyono & Boedi, 2017). Jonassen (2000) stated that learning is a highly complex process. Therefore, an effective learning process cannot be done with a single approach. Upon this ground, ABL was developed as a learning approach that encourages and develops students' active participation in learning theories or concepts through various activities and hands-on experiences in a multitude of learning environments inside and outside of campus. Activity-Based Learning is a procedure where students actively engage in the lesson rather than just sitting, listening, and absorbing the lesson. This experiential learning focuses on learning by doing. Çelik (2018) has it that activity-based learning refers to learning where students physically and mentally explore a subject by simulation of the work environment, manipulation of tools and materials associated with the world of work, or performance of a real work task. Azuka et al. (2013) asserted that an activity-based method is a method of teaching that enables students to be involved in reading, writing, discussion, practical activities, analysis, and evaluation of the topic under discussion. The National Council of Educational Research and Training (NCERT) (2005) cited in Pokhrel (2018) emphasizes that mathematics learning should be facilitated through activities from the very beginning of school education. These activities may involve the use of concrete materials, models, patterns, charts, pictures, posters, games, puzzles, Mathematics laboratory activities, exhibition, and project-based learning and experiments.

A lecture method is a teacher-centered approach whereby the teacher takes part as an active participant and students are at the receiving end most of the time. Fawad (2015) asserted that the lecture method is not only used for teaching theoretical concepts but it is also helpful for giving training of complex skills and procedures. In this method the teacher gives out all the facts he wants the students to know and master, caring very little if at all whether or not, the students are actively participating and contributing to the success of the lesson (Akem, 2007). This method is good for a large class since much work could be easily covered in a shorter time. According to Zohrabi, Torabi, and Baybourdiani, (2012), the main focus of the lecture method is to get students to perform well in promotion examinations rather than other motivation to learning such as innovations and skills to think. Teachers are more dominant sources of information, for instance, all questions posed by students, if any, are answered directly by teachers without students' involvement. In terms of class activities, teachers control all learning experiences. All in all, in the lecture method, the primary task of a teacher is to transfer knowledge to learners.

Noreen and Rana (2019) in their study of the effect of activity-based teaching and lecture method of teaching on students' achievement in the subject of Mathematics found that students taught through activity-based teaching performed better in the post-test. Emaikwu (2012) assessed the relative effectiveness of activity method, discussion method, and lecture methods on students' achievement in secondary school Mathematics and found that activity method performed better than those taught using discussion and lecture methods. Ayhan (2011) investigated the effects of activity-based learning on the success of 8th-grade secondary school students. It was reported that they were able to reach the relationships between given models or data more easily and correctly as the activities progressed. As activities progressed, the students expressed the relationships obtained from the models or problems given in the activities quickly and correctly. In some studies conducted with students, it was found that activity-based teaching increased the students' success in comparison to traditional teaching (Camci, 2012; Küpçü, 2012), it increased students' skills of interpretation better (Camci, 2012). In a similar study, Gürbüz, et al., (2010) compared activity-based learning and traditional learning in terms of their effects on students' development regarding some concepts in the subject of probability. In the experimental study, it was determined that activity-based learning affected the teaching of probability concepts positively. This learning approach in question not only made the process fun but also made learning meaningful. Moreover, Mokiwa and Agbenyeku (2019) investigated the impact of activity-based teaching strategy on gifted students' academic performance and found significant differences in the academic performance between the experimental group exposed to the activity-based teaching strategy when compared to the control group exposed to the traditional lecture method.

Mathematics teachers should employ Bloom's Taxonomy of educational objectives during mathematics teaching and learning. In the 1950s, Benjamin Bloom developed his taxonomy of Educational Objectives. He proposed that learning fitted into one of three psychological domains: Cognitive (which means processing information); Affective (which means attitudes and feelings), and Psychomotor (manipulative or physical) skill. Bloom's Taxonomy, which Benjamin Bloom is best known for, looks at the cognitive domain. This domain categorizes and orders thinking skills and objectives. Bloom's taxonomy is useful to teachers in the sense that it gives teachers a way to think about their teaching and the subsequent learning of their students. Bloom's taxonomy is used to create assessments, evaluate the complexity of assignments, increase the rigor of a lesson, and simplify an activity to help personalize learning, design a summative assessment, plan project-based learning, frame a group discussion, and more. The six levels in Bloom's Taxonomy, according to Heick (2020) are: (1) Remembering level (remember mathematics formulas); (2) Understanding (illustrate the difference between a rectangle and square); (3) Application level; (4) Analysis level (explain how the steps of the scientific process/ mathematical theory work together); (5) Evaluation level (interpret the significance of a given law of Mathematics); (6) Creation level (design a new solution to an 'old' problem that honors/acknowledges the previous failures).

## METHODOLOGY

The study utilized a randomized quasi-experimental design. Ninety-six Senior Secondary II (SSS 2) students for the 2019/2020 academic session in 2 government-owned schools in Calabar Metropolis, Cross River State of Nigeria participated in the study. The students and schools were selected through multistage sampling. The students in their respective schools were randomly assigned to an experimental group and control group respectively; a total of 48 students to activity-based teaching method group and 48 students to lecture teaching method group, making a total of 96 students used in the study.

The instrument for data collection was a theory test (essay-test) consisting of a 20- item Mathematics Test (MT) developed by the researcher and with its items selected from circle geometry and construction which are in the Senior Secondary School II (SSS 2) Mathematics curriculum and the lesson plan were used as instructional tools for the study. Sample of the theory test (essay test) are as follows: (1) A circle of diameter 30cm is cut out from a sheet of square paper of sides 30cm. What is the area of the remaining paper?; (2) Find the perimeter of a circle whose radius is 7/2cm. (3) Construct a right-angled triangle whose base is 12cm and the sum of its hypotenuse and other sides is 18cm; (4) Construct a triangle PQR with its perimeter= 11cm and base angles of $75^{\circ}$ and $30^{\circ}$. The same MT was used as a pretest to determine the readiness levels of the students before experimental implementations, and as a posttest to measure their academic success after the implementations. Face and content validities of the instrument were ensured by two experts in Mathematics Education and one expert in Educational Measurement and Evaluation. The items were spread from comprehension to evaluation levels of Bloom's levels of the cognitive domain. The MT was scored out of 100 marks (5 marks for each item). The reliability coefficient of the instrument for the study was determined by using the Cronbach alpha coefficient and 0.86 was obtained.

The data collected were analyzed in line with the research question and hypothesis. Mean and standard deviation (SD) were used to answer the research question. In a design with pretest, posttest, and a control group, if it is wanted to assess whether the experimental operation was effective or not, the most suitable statistical operation is the method of Analysis of Covariance (ANCOVA), where the pretest is controlled as the mutual variable (Büyüköztürk, 2016). ANCOVA provides a stronger insight into the relationship between the independent and dependent variables by keeping the error variance minimal (Tabachnick & Fidell, 2015). Therefore, Analysis of Covariance (ANCOVA) was used to explore the effect of the teaching methods at p < 0.05 level of significance. SPSS version 22 package software was utilized in the analysis of the data.

## RESULTS

Table 1 shows the performance of students in Mathematics before the treatment (pretest) and after the treatment (posttest). The table shows that the mean score of SSS 2 students before the treatment for the experimental group is 27.19 while their counterparts in the control group are 27.39. Table 2 also shows that after the treatment, the mean score of students in the experimental group increased greatly to 51.71 which is a remarkable improvement while their counterparts in the control group had a low mean score of 29.56 which was very poor. The mathematics performance of the SSS 2 students in the control group was very poor in the pretest and posttest. The results in table 1 show that the activity-based method enhanced the performance of the students while the lecture method did not.

Table 1. The mean pretest and posttest Mathematics scores, and standard deviations scores of the groups

| Group | N | Pretest | | Posttest | |
|---|---|---|---|---|---|
| | | Mean | Standard Deviation | Mean | Standard Deviation |
| Experimental Group (activity-based method) | 48 | 27.19 | 9.56 | 51.71 | 8.84 |
| Control Group (lecture method) | 48 | 27.39 | 9.73 | 29.56 | 6.22 |

ANCOVA was used to determine whether the difference between the posttest Mathematics performance mean scores of students in the experimental (activity-based method) and control (lecture method) groups when the moderating effect of the pretest is controlled was significant. Table 2 shows the results obtained. Table 2 shows that there was a significant difference in the mathematics performance posttest mean scores of the experimental and control groups after controlling for the effect of the pre-test on mathematics scores ($F_{1,114}$, p<0.05.) In other words, after controlling for pretest scores, the posttest success scores are related to the group variable. ANCOVA was used to control the effect of the pretest since the groups were intact classes. Additionally, partial eta squared values provided information on the effects of the Covariant (here, pretest) and the groups on the posttest. The effect size of the group on the posttest was $\eta^2$=0.261, meaning, the group variable explained 26.1% of the posttest variance. On the other hand, the effect size of the covariant (pretest) was even higher as $\eta^2=0.552$, which means, 55.2% of the posttest variance was explained by a pretest.

Table 2. ANCOVA results of difference in posttest mathematics performance mean scores of students in the experimental and control groups when pretest scores were controlled

| Source | Type III Sum of Squares | df | Mean Square | F | p-value | Partial Eta Squared($\eta^2$) |
|---|---|---|---|---|---|---|
| Corrected Model | 9103.117 | 2 | 4551.558 | 73.106 | 0.000 | 0.611 |
| Intercept | 9811.242 | 1 | 9811.242 | 157.587 | 0.000 | 0.629 |
| Pre Maths Score | 2047.607 | 1 | 2047.607 | 32.888 | 0.000 | 0.552 |
| Teaching Methods (post-test) | 7126.293 | 1 | 7126.293 | 114.461 | **0.000** | 0.261 |
| Error | 5790.123 | 93 | 62.259 | | | |
| Total | 154507.000 | 96 | | | | |
| Corrected Total | 14893.240 | 95 | | | | |

## DISCUSSION OF FINDINGS

The finding of this present study revealed that there was a significant difference in the Mathematics performance posttest mean scores of the students taught with the activity-based method (experimental) and students taught with lecture method (control) groups after controlling for the effect of the pre-test on mathematics scores, and a noticeable difference was found to be in favor of the experiment group students in terms of the posttest math academic success scores. In this study, it was found that the experiment group students were academically more successful in comparison to the control group students. Therefore, the activity-based teaching method in Mathematics teaching may increase the academic success of the students in comparison to the lecture method. The reason for this may be that activities in Mathematics classes provide opportunities to work with tangible materials and increase the motivation and interest of the students. Additionally, it may be considered that the presentation and content of the activities made it easier to learn the subject and relate it to daily life. The finding of this study is in agreement with the finding of Mokiwa and Agbenyeku (2019), who found a significant difference in the academic performance between the experimental groups exposed to the activity-based teaching strategy when compared to the control group exposed to the traditional lecture method.

## CONCLUSIONS AND RECOMMENDATIONS

Based on the findings of the present study, it is concluded that activity-based learning increased academic success. Considering the results of relevant studies, it may be considered that activity-based learning will be an applicable approach in all class levels in secondary schools for increasing academic success.

The result of this study has provided an empirical basis that the activity method is an appropriate teaching strategy capable of improving the present dismal achievement of students in Mathematics. This study also contributes to the existing literature on how to improve classroom instructional delivery in Mathematics subject. This study implies that students will benefit greatly from teaching outcomes in schools if the activity method is used as a pedagogical approach in Mathematics instructional delivery. Moreover, the failure rate in Mathematics will decline greatly if an activity-based teaching method is adopted by Mathematics teachers.

It was recommended that secondary school Mathematics teachers should be trained and re-trained to update their knowledge in the use of activity-based teaching for making the teaching and learning of Mathematics more interesting and rewarding. Moreover, Mathematics kits containing materials for activities may be provided to Mathematics teachers to teach Mathematics with activities. A separate room for Mathematics may be allocated in every school and may be equipped with material for activities. Ministries of education, curriculum planners, and developers should outline appropriate activity-based instructional methods for use by teachers to teach any topic highlighted in the Mathematics curriculum. It is recommended for further researchers that studies may be carried out separately for other subjects, especially in science subjects, to know the effectiveness of activity-based teaching.